\documentclass[11pt,a4paper]{amsart}

\newtheorem{theorem}{Theorem}[section]
\newtheorem{lemma}[theorem]{Lemma}
\newtheorem{proposition}[theorem]{Proposition}
\newtheorem{corollary}[theorem]{Corollary}

\newtheorem{question}[theorem]{Question}

\newcommand{\BN}{\mathbb{N}}

\def\inv{^{-1}}
\def\derive{\Rightarrow^{\!{*}}}
\begin{document}

\title{On the rational subset problem for groups}

    \author{MARK KAMBITES}

    \address{Fachbereich Mathematik / Informatik, Universit\"at Kassel,
             34109 Kassel,  Germany
    }
    \email{kambites@theory.informatik.uni-kassel.de}

    \author{PEDRO V.~SILVA}

    \address{Departamento de Matem\'atica Pura,
             Faculdade de Ci\^encias da Universidade do Porto,
             Rua do Campo Alegre 687
             4169-007 Porto, Portugal
     }
    \email{\texttt{pvsilva@fc.up.pt}}

    \author{BENJAMIN STEINBERG}

    \address{School of Mathematics and Statistics, Carleton University,
             Ottawa, Ontario, K1S 5B6, Canada
    }
    \email{bsteinbg@math.carleton.ca}

\begin{abstract}
We use language theory to study the rational subset problem for
groups and monoids. We show that the decidability of this problem is
preserved under graph of groups constructions with finite edge
groups. In particular, it passes through free products amalgamated
over finite subgroups and HNN extensions with finite associated
subgroups. We provide a simple proof of a result of Grunschlag
showing that the decidability of this problem is a virtual property.
We prove further that the problem is decidable for a direct product
of a group $G$ with a monoid $M$ if and only if membership is
uniformly decidable for $G$-automata subsets of $M$. It follows that
a direct product of a free group with any abelian group or
commutative monoid has decidable rational subset membership.
\end{abstract}
\maketitle

\section{Introduction}

A key aspect of combinatorial group theory is the study of
algorithmic \textit{decision problems} in finitely generated groups.
This trend, which was initiated by Dehn \cite{Dehn12}, has always remained
central to the subject, and has recently been given new impetus by interest
in the potential use of group theory as a basis for the development of
secure cryptographic systems (see, for example, \cite{Rivest04,Shpilrain04}).
The best known example of a group-theoretic decision problem is the
\textit{word problem} for a given
finitely generated group $G$: given two words over the generators,
decide if they represent the same element, or equivalently, given a
word over the generators, decide if it represents the identity.

A natural generalization is the \textit{subgroup membership problem}
or \textit{generalized word problem}: given elements $u_1, \dots,
u_n$ and $g$ of $G$ (specified as words over the generators), decide
if $g$ lies in the subgroup generated by $u_1, \dots, u_n$
\cite{Lyndon77}. More generally still, one can consider the
corresponding problem for finitely generated \textit{submonoids} and
\textit{subsemigroups} of groups; Margolis, Meakin and
\v{S}uni\'k~\cite{Margolis05} have recently solved the submonoid
membership problem for a large class of groups by studying
\textit{distortion functions}. Decidability of the subsemigroup
problem is of particular interest since it also implies solvability
of the \textit{order problem}.

From the point of view of computer science and formal language theory,
finitely generated subgroups, submonoids and subsemigroups
are examples of \textit{rational subsets} and it is natural also
to consider the harder \textit{rational subset problem} for $G$: given a
rational subset of $G$ (specified using a finite automaton over the
generators) and a word representing an element of $G$, decide if the latter
belongs to the former.

It is well-known that the subgroup membership problem, and more generally
the rational subset problem, is decidable for free groups and for free
abelian groups \cite{Benois69,Grunschlag99,Lyndon77}.
Many important groups can be built up from these groups using
constructions such as \textit{direct products}, \textit{free products}
(with and without amalgamation), \textit{HNN extensions}
and \textit{graphs of groups} \cite{Lyndon77,Serre03}.
Hence, a natural next step is to consider the extent to which
decidability of these algorithmic problems is preserved under
such operations.

Kapovich, Weidmann and Myasnikov \cite{Kapovich05b} have recently
studied the subgroup membership problem in fundamental groups of
graphs of groups. It is natural to ask whether their results can be
extended to the more general rational subset problem. Their approach
generalizes the well-known \textit{folding} technique of Stallings
\cite{Stallings83}, which starts with a graph easily constructed
from the subgroup generators, and incrementally computes the
important part of the Schreier graph associated to a finitely
generated subgroup of the free group. This method is essentially
automata-theoretic, since the graphs can be viewed as
automata recognising progressively larger parts of the membership
language of the subgroup
\cite{Birget00,Margolis92,Steinberg02}. However, the folding technique relies upon a
degree of symmetry in the automata (the \textit{dual automaton} property,
in the language of \cite{Steinberg02}) which is particular
to the subgroup case, and so does not easily generalize.

The main aim of the present paper is to make a start upon the study
of the relationship between the rational subset problem and the
constructions described above. Our approach is largely
language-theoretic; in Section~\ref{sec_rid}, we introduce an
abstract property of languages which corresponds naturally to
decidability of rational subset membership.  This simple observation
establishes a surprisingly deep connection between algorithmic group
theory and formal language theory. It immediately yields simple
language-theoretic proofs of several non-trivial known results,
including the fact that decidability of the rational subset problem
is a virtual property, which was first proved by Grunschlag
\cite{Grunschlag99}; recall that a virtual property of a group is a
property inherited from finite index subgroups and overgroups. Conversely,
we are also able to use undecidability of a group-theoretic problem to
establish a new undecidability result in language theory: there
exists a context-free language for which it is undecidable which
regular languages it contains.

In Section~\ref{sec_rewrite} we proceed to show that our abstract
property is preserved when taking the ancestor set of a language
under certain infinite rewriting systems. In
Section~\ref{sec_amalg}, we apply these results to show that
decidability of rational subset membership is preserved under graph
of groups constructions \cite{Serre03} with finite edge groups. It
follows in particular that it passes through amalgamated free
products over finite subgroups, and HNN extensions with finite
associated subgroups.

Section~\ref{sec_gauto} introduces a connection between decidability
of rational subsets and the theory of $G$-automata (see for example
\cite{Gilman96,Kambites06}), that is, rational transductions of
group word problems. Specifically, we show that a direct product of
the form $G \times M$ has decidable rational subset membership
exactly if membership is (uniformly) decidable for $G$-automata
subsets of $M$. This combines with a group-theoretic interpretation
\cite{Kambites06} of
a theorem of Chomsky and
Sch\"utzenberger \cite{Chomsky63} and some classical results on
commutative monoids \cite{Eilenberg69,Parikh66} to show that any
direct product of a free group with an abelian group (or commutative
monoid) has decidable rational subset membership.

Finally, in Section~\ref{sec_graph} we consider the subgroup membership
problem and rational subset problem in the important class of
\textit{graph groups} (which are
also known as \textit{right-angled Artin groups}, \textit{trace
  groups}  or \textit{free partially
commutative groups}). We note some consequences for certain of these
groups of our
results from Sections~\ref{sec_amalg} and \ref{sec_gauto}, and pose a number
of questions regarding other graph groups.

While the primary focus of this paper is on groups, the rational subset problem
and associated decision problems are also of interest in more general monoids
and semigroups. For example, recent research of Ivanov, Margolis and Meakin
\cite{Ivanov01} has reduced the word problem for a large class of one-relator
inverse monoids to the submonoid problem in one-relator groups. Since decidability
of rational subset membership is inherited by finitely generated subsemigroups,
many of our results about groups have additional implications for monoids
and semigroups. We also prove directly some new results for monoids,
including the fact that the rational subset problem is decidable for
direct products of free groups with free monoids.

\section{Preliminaries}\label{sec_prelim}

In this section, we provide a brief introduction to rational subsets and
associated decision problems.

\subsection{Rational Subsets, Languages and Transductions}

Let $M$ be a finitely generated monoid. Recall that a \textit{finite
automaton} $P$ over $M$
is a finite directed graph with edges labelled by elements of $M$,
with a distinguished \textit{initial vertex} and a set of distinguished
\textit{terminal vertices}. The labelling of edges extends naturally,
via the multiplication in $M$, to a labelling of paths by elements of
$M$. The \textit{subset recognised} by the automaton is the set of all
elements of $M$ which label paths between the initial vertex and some
terminal vertex. A subset of $M$ recognised by some finite automaton is
called a \textit{rational subset} of $M$.  An alternative description
of the rational subsets of $M$ is as the smallest collection of
subsets of $M$ containing the finite subsets and closed under union,
concatenation and generation of submonoids.

A particularly important case arises when $M = \Sigma^*$ is the free
monoid on an alphabet $\Sigma$, so that the automaton $P$ accepts a
language over $\Sigma$, which we denote $L(P)$. A rational subset of
a free monoid is called a \textit{rational language} or
\textit{regular language}. By \cite[Proposition III.2.2]{Berstel79},
the rational subsets of a monoid $M$ are exactly the homomorphic
images of regular languages. For a detailed introduction to the
theory of regular languages, including a number of alternative
definitions, see \cite{Eilenberg74} or \cite{Hopcroft69}.

Another significant case is when $M = \Sigma^* \times \Omega^*$ is a
direct product of free monoids. A finite automaton over
$\Sigma^* \times \Omega^*$ is called a \textit{finite transducer}
from $\Sigma^*$ to $\Omega^*$. A subset recognised
by a finite transducer, that is, a rational subset of $\Sigma^* \times
\Omega^*$, is called a \textit{rational transduction}.
If $\rho \subseteq \Sigma^* \times \Omega^*$ and $L \subseteq \Sigma^*$
then the \textit{image} of $L$ under $\rho$ is the language
$$L \rho = \lbrace v \in \Omega^* \mid (u, v) \in \rho \text{ for some } u \in L \rbrace \ \subseteq \ \Omega^*.$$
If $\rho$ is a rational transduction then we say that the language $L \rho$
is a \textit{rational transduction of} the
language $L$. Rational transductions are a powerful tool for
studying relationships between formal languages, and many elementary
operations on languages are instances of rational transductions. For
example, given a word $w \in \Sigma^*$ it is straightforward to show
that there is an effectively constructible rational transduction
taking each language $L \subseteq \Sigma^*$ to its \textit{left
translation} by $w$, that is, the language $wL = \lbrace wx \mid x
\in L\rbrace$; an analogous result holds for \textit{right
translation}. More general information about rational transductions
can be found in \cite{Berstel79}.

\subsection{Encodings and Decision Problems}

We shall work with algorithmic problems in which abstract
objects such as monoid elements, automata and languages are regarded as data.
Formally, it is necessary to have an agreed system of
(not necessarily unique) finite encodings for such objects, but for
brevity and clarity it is unhelpful to refer repeatedly to encodings.
Typically we shall introduce an encoding along with the definition of
a class of objects, but thereafter leave the encoding implicit.

All groups and monoids in this paper will be finitely generated, and
we assume that elements are encoded as words over some (except where
otherwise stated, fixed) finite monoid generating set. Our algorithms
are not \textit{uniform} across groups and monoids, and so we do not
need to consider an encoding system for groups and monoids themselves.

Finite automata are presumed to have some suitable encoding, from which
we can extract the vertex set and for any pairs of vertices $p$ and $q$,
the (encodings of) labels of all edges from $p$ to $q$.
Languages and subsets are typically of a type recognised by some kind of
algorithm or automaton or generated by some kind of grammar, and we
assume they are encoded as an automaton, algorithm or grammar.

Now let $\mathcal C$ be a set of subsets of a finitely generated
monoid $M$, and
suppose we have a fixed system of (not necessarily unique) encodings
for elements of $\mathcal C$. For example, $\mathcal C$ might
be the set of rational subsets of $M$, encoded as finite automata. We
say that
\textit{membership is uniformly decidable for $\mathcal C$} if there
is an algorithm which, given an (encoded) set $S \in \mathcal C$, and an
(encoded) element $m \in M$, decides if $m$ lies in $S$.
The following elementary proposition says that this property is
independent of the choice of finite generating set for $M$.
\begin{proposition}\label{prop_genind}
Let $X$ and $Y$ be finite generating sets for a monoid $M$, and let
$\mathcal C$ be a set of subsets of $M$ with a system of encoding for
its elements. Suppose there exists an algorithm which, given an (encoded)
set $S \in \mathcal C$ and a word $u \in X^*$, decides if the word $u$
represents an element of the set $S$. Then there exists an algorithm
which, given an (encoded) set $S \in \mathcal C$ and a word $v \in Y^*$,
decides if the word $v$ represents an element of the set $S$.
\end{proposition}
\begin{proof}
For each letter $y \in Y$, let $w_y \in X^*$ be a word representing
the same element of $M$ as $y$. Define a morphism $\rho : Y^* \to
X^*$ by $y \mapsto w_y$ for all $y \in Y$. Now a word $v \in Y^*$
represents the same element of $M$ as $v \rho \in X^*$. Hence, it
suffices to take the algorithm which, given as input an encoded set
$S \in \mathcal C$ and a word $v \in Y^*$, computes the word $v \rho \in X^*$
by replacing each letter $y$ of $v$ with $w_y$, and then uses the
algorithm given to check whether $v \rho$ represents an element of
$S$.
\end{proof}

One instance of this problem forms the main focus of this paper. A monoid
or group $M$ has \textit{decidable rational subset problem} if membership
is uniformly decidable for the rational subsets of $M$ (encoded as
finite automata over the generating set). Proposition~\ref{prop_genind}
says that this property is invariant under changing the generating set
used to specify the element; much the same argument shows that it is also
invariant under change of the generator set used for the edge labels in
the automaton. Similarly, if a finitely generated monoid $M$ has
decidable rational subset problem, then it is clear that if $N$ is a finitely
generated submonoid of $M$ then the rational subset problem is also decidable
for $N$. Moreover, given an algorithm for the rational subset problem in $M$
and a finite set of generators for the submonoid $N$, we can effectively
compute an algorithm for the rational subset problem in $N$.

Decidability of the rational subset problem for a group $G$ implies
decidability of the word problem and the subgroup, subsemigroup
and submonoid membership problems. Since an element $g \in G$ has
finite order exactly if the subsemigroup generated by $g$ contains
the identity, an algorithm for the rational subset problem also allows
one to decide if a given element has finite order. Moreover, if $g$
does have finite order then, since the word problem is solvable, it
is a simple matter to compute the order by enumerating words representing
powers and checking if they represent the identity.

\subsection{Groups and Word Problems}

If $G$ is a finitely generated group generated by a subset $A$ then the
\textit{word problem} $W_A(G)$ for $G$ with respect to $A$ is the set of
all words in $A^*$ which represent the identity in $G$. More generally,
given a monoid $M$ generated by a subset $A$ and an element $m \in M$,
we denote by $W_A(M, m)$ the set of all words in $A^*$ representing the element $m \in M$.

\section{Regular Intersection Decidability}\label{sec_rid}

In this section we observe that decidability of the rational subset
problem in a group is equivalent to a natural language-theoretic
property of the word problem. By studying the abstract class of
languages with this property, we show that many existing results
about the rational subset problem for groups can be easily deduced
from standard results in language theory. In
Section~\ref{sec_rewrite}, we shall see that this property of
languages is preserved when taking ancestor sets under certain
infinite rewriting systems. In Section~\ref{sec_amalg}, we shall
apply these results to the rational subset problem for fundamental
groups of graphs of groups, and hence for amalgamated free products
and HNN extensions.

We begin by introducing an algorithmic problem that can be associated to
any formal language. Let $L \subseteq \Sigma^*$ be a language. The
\textit{regular intersection decision problem} or \textit{RID problem}
for $L$ is the problem of deciding, given a finite automaton over $\Sigma$,
whether the language recognised intersects (that is, has non-empty
intersection) with $L$. An algorithm
which solves this problem is called an \textit{RID algorithm} for $L$, and
if there exists such an algorithm then $L$ is called \textit{regular
intersection decidable} or \textit{RID}. It is easily seen that an RID
algorithm solves the RID problem for a unique language; hence, RID
languages can be encoded (although not uniquely) as algorithms which
solve their RID problem.

The following provides the connection with the rational subset problem
for groups, and is our motivation for studying RID languages.
\begin{theorem}\label{thm_ridrational}
Let $G$ be a group generated by a finite subset $A$. Then the following
are (effectively) equivalent:
\begin{itemize}
\item[(i)] rational subset membership is decidable for $G$;
\item[(ii)] the word problem for $G$ with respect to $A$ is RID;
\item[(iii)] the RID problem for $W_A(G,g)$ is decidable uniformly in $g \in G$
(where elements of $G$ are encoded as words in $A^*$).
\end{itemize}
\end{theorem}
\begin{proof}
Let $R$ be a rational subset of $G$ encoded as a regular language
$L$ of $A^*$. Then given an element $g\in G$, one has $g\in R$ if
and only if $L$ intersects with $W_A(G,g)$.  This establishes the
equivalence of (i) and (iii).

That (iii) implies (ii) is immediate.  For the converse,
let $w\in A^*$ represent $g\in G$ and let $L\subseteq A^*$ be a
regular language.  Then, since (ii) clearly allows us to solve the
word problem, we can effectively compute a word $u \in A^*$
representing $g\inv$. Moreover, we can effectively compute the
regular language $uL$. Now $L$ intersects with $W_A(G,g)$ if and only
if $W_A(G)$ intersects with $uL$; by assumption, we can test the
latter.
\end{proof}
Despite the simplicity of the proof, we shall see that Theorem~\ref{thm_ridrational}
establishes a surprisingly deep connection between algorithmic group theory
and formal language theory. We observe that the equivalence of (i) and (iii)
holds also for monoids. The following proposition summarises some elementary
properties of the class of RID languages.

\begin{proposition}\label{prop_ridprops}
The class of RID languages is effectively closed under union and
 rational transduction
(and hence also under morphism, inverse morphism, intersection
with regular languages and right and left translation). The
class of RID languages is strictly
contained within the class of recursive languages, and contains the
class of indexed languages (and hence also the context-free and
regular languages).
\end{proposition}
\begin{proof}
If $L$ and $K$ are RID languages and $R$ is a regular language then $R$
intersects with $L \cup K$ exactly if it intersects with $L$ or $K$, so
closure under union is clear.

For closure under rational transduction, suppose $L \subseteq A^*$ is
RID and $\sigma \subseteq A^* \times B^*$ is a rational transduction;
we must show that the language
$$K = L \sigma = \lbrace b \in B^* \mid (a,b) \in \sigma \text{ for some } a \in L \rbrace$$
is RID. To this end, suppose we are given a regular language
$Q \subseteq B^*$. Let
$$P = Q \sigma^{-1} = \lbrace a \in A^* \mid (a,b) \in \sigma \text{ for some } b \in Q \rbrace.$$
Then $P$ is a regular language, and moreover can be effectively computed from
$Q$ \cite[Corollary III.4.2]{Berstel79}. Now it is easily verified that $K$ intersects with $Q$ if and only if
$L$ intersects with $P$; by assumption, the latter can be checked.

The class of indexed languages (encoded as indexed grammars) is
effectively closed under intersection with regular languages
\cite[Corollary~3]{Aho68}, and has decidable emptiness problem
\cite[Theorem~4.1]{Aho68}. Thus, one can test if a regular language
$R$ intersects with an indexed language $L$ by computing the
intersection $L \cap R$, and testing it for emptiness. Hence,
indexed languages are RID. Since regular and context-free
languages are indexed \cite{Aho68}, it follows also that these
languages are RID.

Since singleton sets are (effectively computable as) regular sets, the
membership problem for a language is reducible to the RID problem, so
RID languages are
recursive. On the other hand, the group $F_2 \times F_2$ has solvable word
problem but undecidable rational subset problem
\cite{Lyndon77,Mikhailova58}. Hence, by Theorem~\ref{thm_ridrational},
its word problem is recursive but not RID.
\end{proof}

We now discuss a number of consequences of
Theorem~\ref{thm_ridrational} and Proposition~\ref{prop_ridprops}.
The following lemma will combine with these results to show that the
decidability
of the rational subset problem is inherited by finite index overgroups,
and hence is a \textit{virtual property}. The idea, which
has been used by several authors \cite{Elston04,Gilman96,Holt05}, is
essentially a recoding of the Kaloujnine-Krasner embedding
\cite{Kaloujnine48}. For completeness, we state and prove the lemma explicitly.
\begin{lemma}\label{lem_virtual}
Let $G$ be a finitely generated group and $H$ a finite index
subgroup.  Let $X$ be a finite generating set for $G$ and $Y$ be a
finite generating set for $H$.  Then there is a rational
transduction $\sigma\subseteq Y^*\times X^*$ such that $W_X(G) =
W_Y(H)\sigma$.
\end{lemma}
\begin{proof}
Let
$g_1,\ldots,g_n$ be a complete set of right coset representatives of
$H$ in $G$, assuming without loss of generality that $g_1$ is the
identity.  For each $x\in X$ and $i\in \lbrace 1,\ldots,n\rbrace$,
choose a word $w_{i,x}\in Y^*$ representing the unique element
$h_{i,x}\in H$ such that $g_ix = h_{i,x}g_j$, where $Hg_ix = Hg_j$.
Our transducer has vertex set $G/H$.  The labelled edges are of the
form $Hg_i{\buildrel {(w_{i,x},x)}\over\longrightarrow}Hg_ix$ with
$x\in X$, $i\in \lbrace 1,\ldots,n\rbrace$.  The initial vertex and
terminal vertex are both the coset $H$. If $\sigma\subseteq Y^*\times
X^*$ is the associated rational transduction, then it is easy to see
that $\sigma\inv:X^*\to Y^*$ is a partial function such that, for
$w\in X^*$, $w\sigma\inv$ is defined if and only if $w$ represents
an element of $H$, in which case $w\sigma\inv$ is an element of
$Y^*$ representing the same element of $H$ as $w$; in particular
$w\in W_X(G)$ if and only if $w\sigma\inv \in W_Y(H)$. It follows
that $W_Y(H)\sigma = W_X(G)$.
\end{proof}

Theorem~\ref{thm_ridrational}, Proposition~\ref{prop_ridprops} and
Lemma~\ref{lem_virtual} immediately yield the following result of
Grunschlag~\cite{Grunschlag99}.
\begin{corollary}
Let $G$ and $H$ be finitely generated groups such that $H$ is a
finite index subgroup of $G$ and suppose that $H$ has decidable
rational subset membership. Then $G$ has decidable rational subset
membership.
\end{corollary}

The easier part of a celebrated theorem of Muller and Schupp
\cite{Muller83} states that every finitely generated virtually free
group has context-free word problem. Combining this with
Theorem~\ref{thm_ridrational} and Proposition~\ref{prop_ridprops},
we immediately obtain the following result which subsumes a result
of Benois \cite{Benois69}. A proof of this nature was first
suggested by Margolis and Meakin \cite{MargolisPrivate}.
\begin{corollary}
Finitely generated virtually free groups have decidable rational subset
problem.
\end{corollary}

A significant open question is that of
which finitely generated groups have word problems which are
\textit{indexed languages} \cite{Aho68,Aho69,Gilman98,Lisovik98};
an answer is likely to be of significant interest in both group theory
and language theory.
Theorem~\ref{thm_ridrational} and Proposition~\ref{prop_ridprops} give an
alternative proof of the following result of Lisovik \cite{Lisovik98}, which
gives a necessary condition for a group to have indexed word problem.
\begin{corollary}
Let $G$ be a finitely generated group with indexed word problem. Then $G$
has decidable rational subset membership (and hence decidable subgroup
membership problem and order problem).
\end{corollary}

We can also use Theorem~\ref{thm_ridrational} to obtain some purely
language theoretic results. We have already observed that the group
$F_2 \times F_2$ has undecidable subgroup membership problem
\cite{Lyndon77,Mikhailova58} and hence, by Theorem~\ref{thm_ridrational},
non-RID word problem. In \cite{Holt05}, it is shown that the word problem
of this group is the complement of a context-free language. Combining
these two we obtain the following.
\begin{corollary}\label{cor_cocfnotrid}
The class of complements of context-free languages is not contained in
the class of RID languages.
\end{corollary}

Since a regular language intersects with a language $L$ exactly if
it is not contained in the complement of $L$,
Corollary~\ref{cor_cocfnotrid} is equivalent to the following
statement.
\begin{corollary}
There exists a context-free language $L$ such that there is no algorithm
which decides, given a regular language $R$, whether $R$ is contained in
$L$.
\end{corollary}
We note that the related problems in which the regular language is fixed and
the context-free language varies are also undecidable in general, as a
consequence of undecidability of completeness for context-free languages;
a detailed study of these problems can be found in \cite{Hopcroft69b}.

\section{RID Languages and Monadic Rewriting Systems}\label{sec_rewrite}

In this section, we show that the class of RID languages is closed
under the operation of taking ancestor sets with respect to certain
infinite rewriting systems. In Section~\ref{sec_amalg} we shall use
this result to show that decidability of the rational subset problem
passes through graph of groups constructions with finite edge groups.

A \textit{monadic rewriting system} $\Gamma$ over a finite alphabet $\Sigma$ is a
subset of $\Sigma^*\times (\Sigma\cup \{\varepsilon\})$.  An element
$(w,x)\in \Gamma$ is normally written $w\rightarrow x$.  If
$\mathcal C$ is a class of languages, then $\Gamma$ is called a
\textit{$\mathcal C$-monadic rewriting system} if for each $x\in
\Sigma\cup \{\varepsilon\}$, the set \[\Gamma_x=\{w\in \Sigma^*\mid
(w \rightarrow x)\in \Gamma\}\] belongs to $\mathcal C$. A finite
encoding system for the class $\mathcal C$ naturally gives rise to a
finite encoding system for $\mathcal C$-monadic rewriting systems
which stores for each $x \in \Sigma \cup \lbrace \epsilon \rbrace$
an encoding of $\Gamma_x$. We shall be particularly interested in
the case where $\mathcal C$ is the class of RID languages, encoded
via RID algorithms.

If $\Gamma$ is a monadic rewriting system, then we write
$u\Rightarrow v$ if $u = rws \in \Sigma^*$ and $v=rxs \in \Sigma^*$
with $w\rightarrow x\in \Gamma$.  We denote by $\derive$ the
reflexive, transitive closure of the relation $\Rightarrow$.  If
$u\derive v$, then $v$ is said to be a \emph{descendant}  of $u$
(under $\Gamma$) and $u$ is said be an \emph{ancestor} of $v$ (under
$\Gamma$). If $L\subseteq \Sigma^*$ is a language, then $L\Gamma$
denotes the set of all descendants of $L$ under $\Gamma$ and
$L\Gamma^{-1}$ denotes the set of all ancestors of $L$ under
$\Gamma$.

It is well-known to computer scientists that the set of descendants
of a regular language under a monadic rewriting system is again a
regular language. Moreover, if the rewriting system is finite
\cite{Benois87, Book93} or context-free \cite{Book82} then one can
algorithmically construct an automaton for the language of
descendants. The following theorem gives a more general condition
under which this is possible.

\begin{theorem}\label{thm_monadic}
Let $\Gamma\subseteq \Sigma^*\times (\Sigma\cup\{\varepsilon\})$ be
a monadic rewriting system and let $L\subseteq \Sigma^*$ be a
regular language. Then $L \Gamma$ is regular. Moreover, there is an
algorithm which, given an RID monadic rewriting system (encoded via
the RID algorithms for the $\Gamma_x$, $x\in \Sigma\cup
\{\varepsilon\}$) and a finite automaton recognising a language
$L\subseteq \Sigma^*$, produces an automaton recognising $L\Gamma$.
\end{theorem}
\begin{proof}
Let $M_0$ be a finite automaton recognising a language $L \subseteq
\Sigma^*$. Clearly, by subdividing edges and adding extra vertices
as necessary, we may assume without loss of generality that the
edges of $M_0$ are labelled by elements of $\Sigma \cup \lbrace
\varepsilon \rbrace$. Now starting from $M_0$, we construct a
sequence of automata as follows.

The automaton $M_{i+1}$ has the same vertex set as $M_i$, and all the
edges of $M_i$ (which we call \textit{inherited edges} in
$M_{i+1}$), plus some additional edges (called \textit{new edges})
constructed as follows. For each pair of vertices $p$ and $q$ in $M_i$
denote by $L^i_{pq}$ the set of all words labelling paths from $p$
to $q$ in $M_i$.  For each element $x\in \Sigma\cup \{\varepsilon\}$
such that $\Gamma_x\cap L^i_{pq}\neq \emptyset$, $M_{i+1}$ is given
a new edge from $p$ to $q$ labelled $x$ (unless $M_i$ already had
one).

Moreover, $L^i_{pq}$ is easily seen to be (effectively computable as) a regular language. Hence,
if we are given $\Gamma$ encoded using RID algorithms, we can test whether
each $\Gamma_x$ intersects with $L^i_{pq}$ and so $M_{i+1}$ can be effectively
constructed from $M_i$.

Since every automaton in the sequence has the same vertex set, and at
each stage we only add edges labelled by the (finitely many) letters
in $\Sigma\cup \{\varepsilon\}$ in places where they do not already
exist, the sequence must terminate. That is, there exists $j$ such
that $M_k = M_j$ for all $k \geq j$.  We claim that $L\Gamma =
L(M_k)$.

 Clearly
$L(M_i)\subseteq L(M_{i+1})$ for each $i$ since we have been adding
new edges. We claim that if $v\in L(M_i)$ and $v\Rightarrow w$, then
$w\in L(M_{i+1})$.  Indeed, suppose $v= rus$ and $w=rxs$ with
$u\rightarrow x\in \Gamma$.  Then $v$ labels a successful path $\pi$
in $M_i$ and there is a factorisation $\pi = \rho \Upsilon \sigma$
such that the paths $\rho$, $\Upsilon$, $\sigma$ are respectively labelled by
$r$, $u$, $s$.  Let $e,f$ be the respective initial and terminal
vertices of the path $\Upsilon$.  Then either an edge $\xi$ from $e$
to $f$ labelled by $x$ already exists in $M_i$, or a new edge $\xi$
from $e$ to $f$ with label $x$ is added in the construction of
$M_{i+1}$. Hence $\rho\xi\sigma$ is a successful path labelled by
$w=rxs$ in $M_{i+1}$.  It follows immediately that if $u\in L$ and
$u\derive w$, then $w\in L(M_r)$ where $r$ is the number of steps
needed to derive $w$ from $u$. Hence $L\Gamma\subseteq L(M_k)$.

To show the converse,  it suffices to show that $L(M_{i+1})\subseteq
L(M_i)\Gamma$.  Suppose $w\in L(M_{i+1})$. Then $M_{i+1}$ has a path
$\pi$ from the initial vertex to some terminal vertex labelled $w$;
consider a factorisation $w = w_0 x_1 w_1 x_2 w_2 \dots x_n w_n$
where, in the path $\pi$, each $w_j \in \Sigma^*$ is read along
inherited edges and each $x_j \in \Sigma\cup \{\varepsilon\}$ is
read along a new edge from $e_j$ to $f_j$. Now by the construction
of $M_i$, for each $x_j$ there exists $y_j \in \Sigma^*$ such that
$y_j\rightarrow x_j\in \Gamma$ and $y_j$ labels a path from $e_j$ to
$f_j$ in $M_{i}$. If follows that the word $v=w_0 y_1 w_1 \dots y_n
w_n\in L(M_i)$.  Since $v\derive w$, we see that
$L(M_{i+1})\subseteq L(M_i)\Gamma$. Thus $L(M_k)\subseteq
L(M_0)\Gamma = L\Gamma$, as required.
\end{proof}

In~\cite{Book82}, it is shown that the set of ancestors of a context-free
language under a context-free monadic rewriting system is always
context-free. We obtain an analogous effective result for RID rewriting
systems.
\begin{corollary}\label{cor_ridrewrite}
Let $L\subseteq \Sigma^*$ be an RID language and let $\Gamma$ be an
RID monadic rewriting system over $\Sigma$.  Then $L\Gamma\inv$ is RID.  Moreover,
there is an algorithm which, given an RID language (encoded
as an RID algorithm) and an RID monadic rewriting system $\Gamma$
(encoded as above), outputs an RID algorithm for $L\Gamma\inv$.
\end{corollary}
\begin{proof}
Let $R\subseteq \Sigma^*$.  Then it is straightforward to verify
that $R$ intersects with $L\Gamma\inv$ if and only if $R\Gamma$
intersects with $L$.

The corollary is now immediate from Theorem~\ref{thm_monadic} since
a finite automaton recognising $R\Gamma$ can be effectively
constructed and we can use the RID algorithm for $L$ to check
whether $R\Gamma$ intersects with $L$.
\end{proof}

This corollary allows a new interpretation of the rational subset
problem.
\begin{corollary}\label{cor_ridrewrite2}
Let $G$ be a group with finite generating set $A$.  Then the following are equivalent:
\begin{enumerate}
\item [(i)] $G$ has a decidable rational subset problem;
\item [(ii)] $W_A(G)$ is RID;
\item [(iii)] $W_A(G) = \{\varepsilon\}\Gamma\inv$ for some RID monadic rewriting system $\Gamma$.
\end{enumerate}
\end{corollary}
\begin{proof}
Theorem~\ref{thm_ridrational} gives the equivalence of (i) and (ii).
To show that (ii) implies (iii) simply take the RID monadic rewriting system
consisting of all rules $w\rightarrow \varepsilon$ where $w$ belongs
to the word problem.  The implication (iii) implies (ii) is an
immediate consequence of Corollary~\ref{cor_ridrewrite} and the fact
that singletons are RID languages.
\end{proof}

\section{Graphs of Groups, Amalgamated Products and HNN Extensions}\label{sec_amalg}

In this section, we apply the language-theoretic results of
Section~\ref{sec_rid} to some problems in group theory. We show that
decidability of rational subset membership is preserved under graph
of groups constructions \cite{Serre03} with finite edge groups. A particular
consequence is that this property passes through free products
amalgamated over finite subgroups, and under HNN extensions with
finite associated subgroups.

We briefly recall the definitions of a graph of groups and its
fundamental group; a detailed introduction can be found in
\cite{Serre03}. Let $Y$ be a finite, directed graph with (possibly)
loops and multiple edges. We denote by $V(Y)$ and $E(Y)$ the vertex
and edge sets respectively of $Y$. Let \mbox{$\alpha, \omega : E(Y)
\to V(Y)$} be the functions which take each edge to its start and
end respectively. Suppose we have a fixed-point-free involution $y
\mapsto \overline{y}$ on the edge set $E(Y)$ which is
orientation-reversing, that is, such that $y \alpha = \overline{y}
\omega$ for all $y \in E(Y)$.

A \textit{graph of groups} $(G,Y)$ with underlying graph $Y$ consists of
\begin{itemize}
\item[(i)] for each vertex $v \in V(Y)$, a group $G_v$;
\item[(ii)] for each edge $y \in E(Y)$, a group $G_y$ such that $G_y = G_{\overline{y}}$; and
\item[(iii)] for each edge $y \in E(Y)$, injective morphisms $\alpha_y : G_y \to G_{y \alpha}$ and
$\omega_y : G_y \to G_{y \omega}$ such that $\alpha_y = \omega_{\overline{y}}$ for
all $y \in E(Y)$.
\end{itemize}
We assume that the groups $G_v$ intersect only in the identity, and
that they are disjoint from the edge set $E(Y)$. For each $v \in
V(Y)$, let $\langle X_v \mid R_v \rangle$ be a monoid presentation
for the vertex group $G_v$, with the different generating sets $X_v$ disjoint. Let $B$ denote
the (disjoint) union of $E(Y)$ with all the sets $X_v$. We
define a group $F(G,Y)$ by the monoid presentation
\begin{align*}
F(G,Y) \ = \ \langle \ B \ \mid \ &R_v \ \left(v \in V(Y) \right),\\
&y\overline y  = 1 \ \left(y \in E(Y) \right), \\
&y (g \omega_y) \overline{y} = g \alpha_y \ \left(y \in E(Y), g \in G_y \right) \ \rangle. \\
\end{align*}
Fix a vertex $v_0 \in V(Y)$. We say that a word $w \in B^*$ is of
\textit{cycle type at $v_0$} if it is of the form $w = w_0 y_1 w_1
y_2 w_2 \dots y_n w_n$ where:
\begin{itemize}
\item[(i)] each $y_i \in E(Y)$;
\item[(ii)] $y_1\cdots y_n$ is a path in $Y$ starting and ending at $v_0$;
\item[(iii)] $w_0\in X_{v_0}^*$;
\item[(iv)] for $1\leq i\leq n$, $w_i\in X_{y_i\omega}^*$.
\end{itemize}
The images in $F(G,Y)$ of the words of cycle type at $v_0$ form a
subgroup $\pi_1(G, Y, v_0)$ of $F(G,Y)$, called the
\textit{fundamental group of $(G,Y)$ at $v_0$}.  If the graph $Y$
is connected then the fundamental group is (up to isomorphism)
independent of the choice of vertex $v_0$. The groups
$F(G,Y)$ and $\pi_1(G,Y,v_0)$ are also easily seen to be independent of
the presentations chosen for the vertex groups.

\begin{theorem}\label{thm_gog}
Let $(G,Y)$ be a finite, connected, non-empty graph of finitely
generated groups with finite edge groups. Then the fundamental group
of $(G,Y)$ has decidable rational subset problem if and only if
every vertex group has decidable rational subset problem. Moreover,
the equivalence is effective.
\end{theorem}
\begin{proof}
Since each vertex group embeds into the fundamental group \cite{Serre03},
one implication is immediate.

For the converse, we use the notation defined above. Since the vertex
groups are assumed to be finitely generated, we may assume that the
generating sets $X_v$ are finite. Moreover, since the edge groups are
finite and there are finitely many edges, we may assume without loss of
generality that for every
edge $y$, the sets $X_{y \alpha}$ and $X_{y \omega}$ contain a
letter representing each non-identity element of
$G_{y\alpha}\alpha_y$ and $G_{y\omega}\omega_y$, respectively. Let
$B$ be the (disjoint) union of all the sets $X_v$ and the edges of
$Y$. Then $B$ is a finite generating set for the group $F(G,Y)$.

Fix a vertex $v_0 \in V(Y)$, and let $P \subseteq B^*$ denote the set of
words of cycle type at
$v_0$. We claim first that the intersection of $P$ with the word
problem $W$ of $F(G,Y)$ is RID. To show this we define an RID
monadic rewriting system $\Gamma$ over $B$ with the following three
types of rules:
\begin{itemize}
\item  If $v$ is a vertex and $w\in X_v$ represents the identity of $G_v$, then there is
a rule $w\rightarrow \varepsilon$;
\item If $y$ is an edge, then there is a rule $y\overline{y}\rightarrow
\varepsilon$;
\item  If $y$ is an edge and $h\in X_{y\alpha}$ is a letter representing $g\in
G_y\alpha _y$ and $w\in X_{y\omega}^*$ represents $g\omega_y$, then
there is a rule $yw\overline y\rightarrow h$.
\end{itemize}

Since all the vertex groups $G_v$ are assumed to be RID, it follows easily
from Theorem~\ref{thm_ridrational} and Proposition~\ref{prop_ridprops}
that $\Gamma$ is an RID monadic rewriting system. Let $L = \lbrace
\varepsilon \rbrace \Gamma\inv$. Since a singleton language is RID,
we deduce by Corollary~\ref{cor_ridrewrite} that $L$ is RID (as
singletons are clearly RID languages). We claim that $L$ is the
intersection of the language $P$ of words of cycle type at
$v_0$ with the word problem $W$ of $F(G,Y)$.

Clearly, the rewriting rules in $\Gamma$ are relations satisfied in
$F(G,Y)$, from which it follows that $L \subseteq W$. It is also easy
to see that the rewriting
rules in $\Gamma$ and their inverses preserve paths of cycle type at
$v_0$, so that $L \subseteq P$. Thus, $L \subseteq W \cap P$.

Conversely, suppose that $w\in W \cap P$.  Then we can write \[w =
w_0 y_1 w_1 y_2 w_2 \dots y_n w_n\] as per (i)-(iv) above. We
proceed by induction on the parameter $n$ (which we shall term the
length of $w$) to show that $w\in L$. If $n=0$, then $w$ represents
$1$ in $G_{v_0}$ and so $w\rightarrow \varepsilon \in \Gamma$,
whence $w\Rightarrow \varepsilon$, establishing $w\in L$.

Assume that all elements of $W\cap P$ of length at most $k<n$ belong
to $L$. Since $w\in W\cap P$, it follows by
\cite[Theorem~I.11]{Serre03} that there exists $i$ such that
$y_{i+1} = \overline{y_{i}}$ and $w_i$ represents an element $g
\omega_{y_i}$ for some $g \in G_{y_i}$. There are two cases: either
$g=1$ or $g\neq 1$.  If $g=1$, then we may apply the rule
$w_i\rightarrow \varepsilon$ followed by the rule
$y_i\overline{y_i}\rightarrow \varepsilon$ to show that if $w'=w_0
y_1 \dots w_{i-1}w_{i+1}y_{i+2}\cdots w_n$, then  $w\derive w'$.  As
$w'\in W\cap P$ and has smaller length, it follows by induction that
$w'\in L$, from which we obtain $w\in L$. If $g\neq 1$, let $h\in
X_{y_i\alpha}$ be a letter representing $g\alpha_{y_i}$. Then by
definition $y_i w_i \overline{y_{i}}\rightarrow h$ lies in $\Gamma$.
Since $y_{i+1} = \overline{y_i}$, if we set $w' = w_0 y_1 \dots
w_{i-1} h w_{i+1}y_{i+2}\cdots w_n$, then $w\Rightarrow w'$. Since
$w'\in W\cap P$ represents the identity in $F(G,Y)$ and its
expression is shorter, we obtain by induction that $w'\derive
\varepsilon$  and so $w\derive\varepsilon$.   Thus, we conclude that
$W \cap P=L$ and so is RID, as claimed.

Now let $X$ be a finite monoid generating set for the fundamental
group $\pi_1(G, Y, v_0)$; this group is indeed finitely generated,
namely by a set in correspondence with the disjoint union of the
generating sets $X_v$
and the set of edges not belonging to some spanning tree for
$Y$~\cite{Serre03}. Choose a morphism $\rho : X^* \to B^*$ which
takes each element $x \in X$ to some word $w_x \in P$ which
represents the same element of $\pi_1(G,Y,v_0)$ as $x$. Since $P$ is
a submonoid of $B^*$, the image of $\rho$ lies in $P$. Hence, for
any word $w \in X^*$, $w$ lies in the word problem of
$\pi_1(G,Y,v_0)$ with respect to $X$ if and only if $w \rho$ lies in
$W$, that is, if and only if $w \rho$ lies in $W \cap P = L$. So the
word problem for $\pi_1(G,Y,v_0)$ is an inverse morphic image of an
RID language, and so by Proposition~\ref{prop_ridprops} is RID. The
theorem now follows from Theorem~\ref{thm_ridrational}.
\end{proof}

An \textit{HNN extension} is the fundamental group of a graph of
groups with a single vertex $v$ and a matched pair $y, \overline{y}$
of loops at $v$. The vertex group $G_v$ is the base group of the HNN
construction, while the edge group $G_y = G_{\overline{y}}$ is
isomorphic to the associated subgroups~\cite{Serre03}. Hence, we
obtain the following corollary.
\begin{corollary}\label{cor_hnn}
Let $H$ be an HNN extension of a group $G$ with finite associated subgroups.
Then $H$ has decidable rational subset problem if and only if $G$
has decidable rational subset problem. Moreover, the equivalence is
effective.
\end{corollary}

Similarly, an \textit{amalgamated free product} corresponds to the
fundamental group of a graph of groups with two vertices connected
by a matched pair of edges~\cite{Serre03}, yielding the following
result.
\begin{corollary}\label{cor_amalg}
Let $H$ be a free product of groups $G_1$ and $G_2$ amalgamated over a
finite subgroup. Then $H$ has decidable rational subset problem if
and only if $G_1$ and $G_2$ have decidable rational subset problem.
Moreover, the equivalence is effective.
\end{corollary}

We remark that the proof of Theorem~\ref{thm_gog} establishes a more
general fact. Recall that if $\mathcal C$ is a class of languages closed
under inverse morphisms then the property of a finitely generated group
having word problem in $\mathcal C$ is independent of the finite generating
set chosen \cite[Lemma 1]{Holt05}.

\begin{theorem}\label{thm_goggeneral}
Let $\mathcal C$ be a class of languages closed under inverse
morphism, left and right translation and closed under taking the
ancestors of $\{\varepsilon\}$ for any $\mathcal C$-monadic
rewriting system.  Then the fundamental group of a finite graph of
finitely generated groups with finite edge groups has word problem
in $\mathcal C$ if all the vertex groups have word problems
in $\mathcal C$. The converse is true if $\mathcal C$ is
closed under intersection with regular languages.
\end{theorem}

The assumptions of Theorem~\ref{thm_goggeneral} are satisfied by, for
example, the
class of context-free languages~\cite{Berstel79,Book82}. A
well-known theorem of Muller and Schupp \cite{Muller83} augmented by
a subsequent result of Dunwoody \cite{Dunwoody85} says that a group
has context-free word problem if and only if it is virtually free.
Our method therefore gives the following result, which can also be
proved in many other ways, for instance from work of Karrass,
Pietrowski and Solitar \cite{Karrass73}.
\begin{corollary}
Let $(G,Y)$ be a finite, connected graph of groups with finite edge
groups. Then the fundamental group of $(G,Y)$ is virtually free if
and only if every vertex group is virtually free.
\end{corollary}

\section{$G$-automata and Rational Subsets}\label{sec_gauto}

In this section, we consider the rational subset problem in direct
products. We show that if $G$ is a finitely generated group and $M$
a finitely generated monoid, then the rational subset problem for $G
\times M$ is decidable exactly if the subsets of $M$ defined by
\textit{$G$-automata} have uniformly decidable membership problem.
This combines with a group-theoretic interpretation \cite{Kambites06} of
a well-known theorem of Chomsky and
Sch\"utzenberger \cite{Chomsky63} to give a characterisation of
rational subset membership in direct products of the form $F \times M$
with $F$ a free group, in terms of the uniform decidability of
membership for context-free subsets of $M$.

Let $G$ be a finitely generated group. Recall that a \textit{$G$-automaton
over the alphabet $\Sigma$} is a finite automaton $P$ over $G\times \Sigma^*$.
The \emph{$G$-automaton language accepted by $P$} is the set of all words
$w\in \Sigma^*$ such that $(1,w)$ belongs to the rational subset recognised
by $P$ \cite{Gilman96,Kambites06}. The $G$-automaton languages are exactly
the rational transductions of the word problem of $G$ \cite[Proposition~2]{Kambites06}.
More generally, we say that a
\textit{$G$-automaton over a monoid $M$} is a finite automaton $P$
over $G \times M$. The \textit{$G$-automaton subset} recognised by $P$ is
the set of all elements $m \in M$ such that $(1,m)$ belongs to the rational
subset recognised by $P$. It is easily seen that the $G$-automaton subsets
are exactly the homomorphic images of $G$-automaton languages.

Having fixed
$G$ and $M$ and some finite generating sets $A$ and $B$ respectively, the
$G$-automata subsets of $M$ have a natural encoding as finite automata
over $G \times M$. Thus, we may ask whether membership is uniformly
decidable for $G$-automaton subsets of $M$; an argument similar to the
proof of Proposition~\ref{prop_genind} shows that this property is
independent of the choice of generating sets.
The following result relates decidability properties of
$G$-automaton subsets to the rational subset membership problem.

\begin{theorem}\label{thm_gautodec}
Let $G$ be a finitely generated group, and $M$ a finitely generated
monoid. Then the following are equivalent:
\begin{itemize}
\item[(i)] the rational subset problem for $G\times M$ is decidable;
\item[(ii)] membership is uniformly decidable for $G$-automaton subsets of $M$.
\end{itemize}
\end{theorem}
\begin{proof}
Let $A$ and $B$ be finite generating sets for $G$ and $M$ respectively, so
that $A \cup B$ is a generating set for $G \times M$. In view of our comments
above, we may assume that all words and automata are encoded using these
generating sets.

Suppose first that (i) holds, that is, that the rational subset problem
is decidable, and that we are given a finite automaton $P$ over
$G \times M$ and an element $m \in M$. In view of our
choice of generators, the word over $B$ encoding $m \in M$ also encodes
the element $(1,m) \in G \times M$. Now $m$ lies in the $G$-automaton
language defined by $P$ if and only if $(1,m)$ lies in the rational
subset defined by $P$. By assumption, this can be tested, so that (ii)
holds

Conversely, suppose (ii) holds, and that we are given a finite automaton
$P$ defining a rational subset of $R \subseteq G \times M$ and an element
$(g,m) \in G \times M$. Since $(g,m)$ is encoded as a word over $A \cup B$,
we can easily compute a word representing $(g,1)$ and from that, a word
representing $(g^{-1},1)$. It follows that we can construct from $P$ a
finite automaton $Q$ recognising the rational subset $(g^{-1},1) R$. Now
$(g,m)$ lies in $R$ if and only if $(1,m)$ lies in $(g^{-1},1) R$, that is,
if and only if $(1,m)$ is accepted by $Q$ as a $G$-automaton. Once again,
this can be tested, which shows that (i) holds and completes the proof.
\end{proof}

We note that, since property (i) in the statement of
Theorem~\ref{thm_gautodec} is symmetric in $G$ and $M$, we obtain
the following corollary for $G$-automaton subsets of groups.
\begin{corollary}
Let $G$ and $H$ be finitely generated groups. Then the following are
equivalent:
\begin{itemize}
\item[(i)] membership is uniformly decidable for $G$-automata subsets of $H$;
\item[(ii)] membership is uniformly decidable for $H$-automata subsets of $G$.
\end{itemize}
\end{corollary}

We now turn our attention to the implications of
Theorem~\ref{thm_gautodec} in the case that the group $G$ is a finitely
generated free group $F$ of rank $2$ or more. 
By \cite[Theorem~7]{Kambites06}, which is essentially a group-theoretic
restatement of the Chomsky-Schutzenberger theorem \cite{Berstel79,Chomsky63},
the languages accepted by $F$-automaton are exactly the context-free
languages, Combining with Theorem~\ref{thm_gautodec}, we immediately
obtain the following corollary, where $F_n$ denotes a free group of rank
$n$.
\begin{corollary}\label{thm_cfrat}
Let $M$ be a finitely generated monoid. Then the following are
equivalent:
\begin{itemize}
\item[(i)] the rational subset problem is decidable for the direct
product $F_2\times M$;
\item [(ii)] the rational subset problem is decidable for the direct
product $F_n\times M$ for all $n\geq 0$;
\item[(iii)] membership is uniformly decidable for context-free subsets of $M$
(encoded as context-free grammars over a finite generating set).
\end{itemize}
\end{corollary}

This leads to a new proof of the following result of Frougny,
Sakarovitch and Schupp~\cite{Frougny89}.

\begin{corollary}
The membership problem for context-free subsets of a free
non-abelian group on at least $2$ generators is undecidable.
\end{corollary}
\begin{proof}
If $F$ is a free non-abelian group then $F\times F$ has undecidable
subgroup membership problem~\cite{Lyndon77,Mikhailova58} and hence
undecidable rational subset problem. The result now follows from
Corollary~\ref{thm_cfrat}.
\end{proof}

Applying known results from language theory, we obtain the
following.
\begin{theorem}\label{thm_ratdec}
Let $F$ be a free group.   The rational subset problem is decidable
for:
\begin{enumerate}
\item[(i)] direct products $F\times A$ with $A$ a finitely generated
abelian group;
\item[(ii)] direct products $F\times M$ with $M$ a finitely generated commutative monoid; and
\item[(iii)] direct products $F\times X^*$ with $X$ a finite set.
\end{enumerate}
\end{theorem}
\begin{proof}
In each case, by Corollary~\ref{thm_cfrat}, it suffices to show that
membership is uniformly decidable for context-free
subsets of the monoid in question (encoded as context-free grammars over
a finite generating set).  This is well-known for the case
of a free monoid $X^*$ \cite[Section 6.3]{Hopcroft69}, so (iii) holds.

Since abelian groups are examples of commutative monoids, it
suffices now to prove case (ii). Let $M$ be a commutative monoid
generated by a finite subset $X$; then there is a surjective
morphism $\rho : X^* \to M$. Let $\sigma : X^* \to \BN^X$ be the
canonical morphism from the free monoid $X^*$ to the free
commutative monoid $\BN^X$ on $X$; clearly $\rho$ factors through
$\sigma$ via a morphism $\tau:\BN^X\to M$.

Now suppose we are given a context-free subset of $M$ (encoded as a
context-free grammar over $X$) and an element $m \in M$ encoded as a
word $w \in X^*$. Let $L$ be the language over $X$ generated by the
grammar. Then
by Parikh's Theorem \cite{Parikh66}, we can effectively compute a
regular language $L'$ such that $L'\sigma = L\sigma$. Hence $L\rho =
L'\rho$ is a rational subset of $M$.  But a result of Eilenberg and
Sch\"utzenberger~\cite{Eilenberg69} shows that the preimage of any
rational subset of $M$ under $\tau$ is a rational subset of $\BN^X$.
Moreover, the proof is effective: one can effectively find a regular
language $L''\subseteq X^*$ such that $L''\sigma = L'\rho\tau\inv$.
So $m\in L$ if and only if $w\sigma \in L''\sigma$.  Hence we have
reduced our problem to the rational subset membership problem for
$\BN^X$. Now $\BN^X$ is a finitely generated submonoid of
$\mathbb{Z}^X$.  It was observed by Grunschlag~\cite{Grunschlag99}
that the description of rational subsets in commutative monoids as
semilinear sets, due to Eilenberg and
Sch\"utzenberger~\cite{Eilenberg69}, leads to an immediate
solution of the rational subset membership problem for
$\mathbb{Z}^X$ (and hence $\BN^X$) via integer programming.  This
completes the proof.
\end{proof}

\section{Graph Groups}\label{sec_graph}

In this section, we briefly discuss the subgroup membership problem
and rational subset problem for the class of graph groups (which are
also known as right-angled Artin groups, trace groups or free
partially commutative groups).

Let $\Gamma$ be a finite undirected graph; we denote by $V(\Gamma)$
the vertex set of $\Gamma$, and by $E(\Gamma)$ the edge set of
$\Gamma$, which we view as a symmetric, reflexive subset of
$V(\Gamma) \times V(\Gamma)$. Recall that the graph group
$G(\Gamma)$ is the group with presentation
$$\langle \ V(\Gamma) \ \mid \ ef = fe \text{ for all } (e, f) \in E(\Gamma) \ \rangle.$$

The extreme examples of graph groups, obtained when $\Gamma$ has no edges
or is complete, are free groups and free abelian groups respectively. In
general,
there are many properties of groups which are easily seen to hold for free
groups and free abelian groups, but for radically different reasons.
Establishing the extent to which such properties hold in general graph
groups is often much more difficult, but can be very enlightening.
Decidability of the subgroup membership problem and decidability of the
rational subset problem are two such properties.

A recent result of Kapovich, Weidmann and Myasnikov
\cite{Kapovich05b} shows that the subgroup membership problem is
decidable for graph groups on finite graphs without chord-free
cycles of length four or more. On the other hand, the graph group on a four-cycle is a
direct product of non-abelian free groups; it follows from the
results of Mikhailova \cite{Lyndon77,Mikhailova58} discussed above
that the subgroup membership problem is undecidable for $G(\Gamma)$
whenever $\Gamma$ contains a chord-free four-cycle. Decidability of
the subgroup membership problem seems to be open for graph groups
$G(\Gamma)$ where $\Gamma$ contains chord-free cycles but not of
length four; these groups do not contain a direct product of
non-abelian groups as a subgroup \cite{Kambites06b}, and so
Mikhailova's result does not assist. The following question is an
obvious starting point for research in this direction.
\begin{question}
Is the subgroup membership problem decidable for the graph group on an
$n$-cycle with $n \geq 5$?
\end{question}

Similarly, we have seen that the rational subset problem is decidable
for free groups
and free abelian groups. Combining Corollary~\ref{cor_amalg} and
Theorem~\ref{thm_ratdec} we immediately obtain decidability for a
somewhat larger class.
\begin{corollary}\label{cor_graph}
The rational subset membership problem is decidable for free
products of direct products of a free group with a free abelian
group.
\end{corollary}
In graph-theoretic terms, Corollary~\ref{cor_graph} applies to $G(\Gamma)$
where every connected component of $\Gamma$ is the join
(see \cite{Droms92}) of a complete graph and a graph with no edges. On the other
hand, we saw above that the subgroup membership problem is undecidable for
$G(\Gamma)$ where $\Gamma$ contains a four-cycle without chords,
so the rational subset problem is also undecidable in these cases.
Many cases remain open; these are summarised in the
following question.
\begin{question}
Is the rational subset membership problem decidable for $G(\Gamma)$ when
\begin{itemize}
\item[(i)] $\Gamma$ is a three-edge line?
\item[(ii)] $\Gamma$ has no chord-free cycles of length four or more?
\item[(iii)] $\Gamma$ is an $n$-cycle with $n \geq 5$?
\item[(iv)] $\Gamma$ contains no chord-free four-cycles?
\end{itemize}
\end{question}

\section*{Acknowledgements}

The authors would like to thank Stuart Margolis and Friedrich Otto
for some helpful conversations. The research of the first author was
supported by a
Marie Curie Intra-European Fellowship within the 6th European
Community Framework Programme; he would also like to thank Kirsty
for all her support and encouragement. The second author
acknowledges support from C.M.U.P., financed by F.C.T. (Portugal)
through the programmes POCTI and POSI,  with national and European
Community structural funds, as well as the European Science
Foundation programme AutoMathA.  The work of the third author was
supported by an NSERC discovery grant.

\bibliographystyle{plain}

\def\cprime{$'$} \def\cprime{$'$}

\end{document}